\documentstyle[twoside]{article}
\include{graphicx}
\input epsf
\title{Some remarks on the theorems of Wright and Braaksma on the Wright function ${}_p\Psi_q(z)$
\footnote{This paper is partly based on the internal report \cite{P09}.}}
\author{\sc R. B.\ Paris \\
{\em University of Abertay Dundee, Dundee DD1 1HG, UK}}
\begin{document}
\def\f#1#2{\mbox{${\textstyle \frac{#1}{#2}}$}}
\def\dfrac#1#2{\displaystyle{\frac{#1}{#2}}}
\def\boldal{\mbox{\boldmath $\alpha$}}
\newcommand{\bee}{\begin{equation}}
\newcommand{\ee}{\end{equation}}
\newcommand{\lam}{\lambda}
\newcommand{\ka}{\kappa}
\newcommand{\al}{\alpha}
\newcommand{\fr}{\frac{1}{2}}
\newcommand{\fs}{\f{1}{2}}
\newcommand{\g}{\Gamma}
\newcommand{\br}{\biggr}
\newcommand{\bl}{\biggl}
\newcommand{\ra}{\rightarrow}
\newcommand{\mbint}{\frac{1}{2\pi i}\int_{c-\infty i}^{c+\infty i}}
\newcommand{\mbcint}{\frac{1}{2\pi i}\int_C}
\newcommand{\mboint}{\frac{1}{2\pi i}\int_{-\infty i}^{\infty i}}
\newcommand{\gtwid}{\raisebox{-.8ex}{\mbox{$\stackrel{\textstyle >}{\sim}$}}}
\newcommand{\ltwid}{\raisebox{-.8ex}{\mbox{$\stackrel{\textstyle <}{\sim}$}}}
\renewcommand{\topfraction}{0.9}
\renewcommand{\bottomfraction}{0.9}
\renewcommand{\textfraction}{0.05}
\newcommand{\mcol}{\multicolumn}
\date{}
\maketitle
\pagestyle{myheadings}
\markboth{\hfill \sc R. B.\ Paris  \hfill}
{\hfill \sc Remarks on the theorems of Wright and Braaksma \hfill}
\begin{abstract}
We carry out a numerical investigation of the asymptotic expansion of the so-called Wright function ${}_p\Psi_q(z)$
(a generalised hypergeometric function) in the case when exponentially small terms are present. This situation is covered by two theorems of Wright and Braaksma. We demonstrate that a more precise understanding of the behaviour of ${}_p\Psi_q(z)$ is obtained by taking into account the Stokes phenomenon.
\end{abstract}
\vspace{0.3cm}

\noindent $\,$\hrulefill $\,$

\vspace{0.2cm}

\begin{center}
{\bf 1. \  Introduction}
\end{center}
\setcounter{section}{1}
\setcounter{equation}{0}
\renewcommand{\theequation}{\arabic{section}.\arabic{equation}}
We consider the Wright function (a generalised hypergeometric function) defined by
\bee\label{e11}
{}_p\Psi_q(z)\equiv{}_p\Psi_q\bl(\!\!\begin{array}{c}(\alpha_1,a_1), \ldots ,(\alpha_p,a_p)\\(\beta_1, b_1), \ldots ,(\beta_q,b_q)\end{array}\!\!;z\!\br)=\sum_{n=0}^\infty g(n)\,\frac{z^n}{n!}, 
\ee
\bee\label{e11a}
g(n)=\frac{\prod_{r=1}^p\Gamma(\alpha_rn+a_r)}{\prod_{r=1}^q\Gamma(\beta_rn+b_r)},
\ee
where $p$ and $q$ are nonnegative integers, the parameters $\alpha_r$  and 
$\beta_r$ are real and positive and $a_r$ and $b_r$ are
arbitrary complex numbers. We also assume that the $\alpha_r$ and $a_r$ are subject to 
the restriction
\bee\label{e11ab}
\alpha_rn+a_r\neq 0, -1, -2, \ldots \qquad (n=0, 1, 2, \ldots\ ;\, 1\leq r \leq p)
\ee
so that no gamma function in the numerator in (\ref{e11}) is singular.
In the special case $\alpha_r=\beta_r=1$, the function ${}_p\Psi_q(z)$ reduces to a multiple of 
the ordinary hypergeometric function 
\[{}_p\Psi_q(z)=\frac{\prod_{r=1}^p\g(a_r)}{\prod_{r=1}^q\g(b_r)}\,{}_pF_q\bl(\begin{array}{c}a_1, \ldots ,a_p\\b_1, \ldots ,b_q\end{array}\!\!;z\br);\]
see, for example, \cite [ p.~40]{S}.

We introduce the parameters associated\footnote{Empty sums and products are to be interpreted as zero and unity, respectively.} with $g(n)$ given by
\[\kappa=1+\sum_{r=1}^q\beta_r-\sum_{r=1}^p\alpha_r, \qquad 
h=\prod_{r=1}^p\alpha_r^{\alpha_r}\prod_{r=1}^q\beta_r^{-\beta_r},\]
\bee\label{e12}
\vartheta=\sum_{r=1}^pa_r-\sum_{r=1}^qb_r+\f{1}{2}(q-p),\qquad \vartheta'=1-\vartheta.
\ee
If it is supposed that $\alpha_r$ and $\beta_r$ are such that $\kappa>0$ then ${}_p\Psi_q(z)$ 
is uniformly and absolutely convergent for all finite $z$. If $\kappa=0$, the sum in (\ref{e11})
has a finite radius of convergence equal to $h^{-1}$, whereas for $\kappa<0$ the sum is divergent 
for all nonzero values of $z$. The parameter $\kappa$ will be found to play a critical role 
in the asymptotic theory of ${}_p\Psi_q(z)$ by determining the sectors in the $z$-plane 
in which its behaviour is either exponentially large, algebraic or exponentially small 
in character as $|z|\ra\infty$.

The determination of the asymptotic expansion of ${}_p\Psi_q(z)$ for $|z|\ra\infty$ and finite 
values of the parameters has a long history; for details, see \cite[ \S 2.3]{PW}.
Detailed investigations were carried out by Wright \cite{W1, W2} and by
Braaksma \cite{Br} for a more general class of integral functions than (\ref{e11}). We present a summary of their results related to the asymptotic expansion of ${}_p\Psi_q(z)$ for large $|z|$ in Section 2. Our purpose here is to consider two of the expansion theorems involving the presence of exponentially small expansions valid in certain sectors of the $z$-plane. We demonstrate by numerical computation that a more precise understanding of the asymptotic structure of ${}_p\Psi_q(z)$ can be achieved by taking into account the Stokes phenomenon.

\vspace{0.6cm}

\begin{center}
{\bf 2. \ Standard asymptotic theory for $|z|\ra\infty$}
\end{center}
\setcounter{section}{2}
\setcounter{equation}{0}
\renewcommand{\theequation}{\arabic{section}.\arabic{equation}}
We first state the standard asymptotic expansion of the integral function ${}_p\Psi_q(z)$ as 
$|z|\ra\infty$ for $\kappa>0$ and finite values of the parameters given in \cite{W2} and 
\cite{Br}; see also \cite[ \S2.3]{PK}.
To present this expansion we introduce the exponential expansion $E_{p,q}(z)$ and the 
algebraic expansion $H_{p,q}(z)$ associated with ${}_p\Psi_q(z)$.

The exponential expansion $E_{p,q}(z)$ can be obtained from the Ford-Newsom theorem \cite{F,N}. A simpler derivation of this result in the case ${}_p\Psi_q(z)$  based on the Abel-Plana form of the well-known Euler-Maclaurin summation formula is given in \cite[pp.~42--50]{PW}.
We have the formal asymptotic sum
\bee\label{e22c}
E_{p,q}(z):=Z^\vartheta e^Z\sum_{j=0}^\infty A_jZ^{-j}, \qquad Z=\kappa (hz)^{1/\kappa},
\ee
where the coefficients $A_j$ are those appearing in the inverse factorial expansion of $g(s)/s!$ given by  
\bee\label{e22a}
\frac{g(s)}{\g(1+s)}=\kappa (h\kappa^\kappa)^s\bl\{\sum_{j=0}^{M-1}\frac{A_j}{\Gamma(\kappa s+\vartheta'+j)}
+\frac{\rho_M(s)}{\Gamma(\kappa s+\vartheta'+M)}\br\}.
\ee
Here $g(s)$ is defined in (\ref{e11a}) with $n$ replaced by $s$, $M$ is a positive integer and $\rho_M(s)=O(1)$ for $|s|\ra\infty$ in $|\arg\,s|<\pi$.
The leading coefficient $A_0$ is specified by
\bee\label{e22b}
A_0=(2\pi)^{\frac{1}{2}(p-q)}\kappa^{-\frac{1}{2}-\vartheta}\prod_{r=1}^p
\alpha_r^{a_r-\frac{1}{2}}\prod_{r=1}^q\beta_r^{\frac{1}{2}-b_r}.
\ee
The coefficients $A_j$ are independent of $s$ and depend only on the parameters $p$, $q$, $\alpha_r$, 
$\beta_r$, $a_r$ and $b_r$. An algorithm for their evaluation is described in the appendix.

The algebraic expansion $H_{p,q}(z)$ follows from the Mellin-Barnes integral representation \cite[ \S 2.4]{PK}
\bee\label{e24aa}
{}_p\Psi_q(z)=\frac{1}{2\pi i}\int_{-\infty i}^{\infty i} \Gamma(s)g(-s)(ze^{\mp\pi i})^{-s}ds,\qquad |\arg(-z)|<\pi(1-\fs\kappa),
\ee
where the path of integration is indented near $s=0$ to separate\footnote{This is always 
possible when the condition (\ref{e11ab}) is satisfied.} the poles of $\g(s)$ from those of 
$g(-s)$ situated at 
\bee\label{e24a}
s=(a_r+k)/\alpha_r, \qquad k=0, 1, 2, \dots\, \ (1\leq r\leq p).
\ee
In general there will be $p$ such sequences of simple poles though, depending on the values 
of $\alpha_r$ and $a_r$, some of these poles could be multiple poles or even ordinary 
points if any of the $\Gamma(\beta_rs+b_r)$ are singular there. Displacement of the contour to the 
right over the poles of $g(-s)$ then yields the algebraic expansion of 
${}_p\Psi_q(z)$ valid in the sector in (\ref{e24aa}). 

If it is assumed that the parameters are such that 
the poles in (\ref{e24a}) are all simple we obtain the algebraic expansion given by 
$H_{p,q}(z)$, where
\bee\label{e25}
H_{p,q}(z):=\sum_{m=1}^p\alpha_m^{-1}z^{-a_m/\alpha_m}S_{p,q}(z;m)
\ee
and $S_{p,q}(z;m)$ denotes the formal asymptotic sum
\bee\label{e25a}
S_{p,q}(z;m):=\sum_{k=0}^\infty \frac{(-)^k}{k!}\Gamma\left(\frac{k+a_m}{\alpha_m}\right)\,
\frac{\prod_{r=1}^{'\,p}\Gamma(a_r-\alpha_r(k+a_m)/\alpha_m)}
{\prod_{r=1}^q\Gamma(b_r-\beta_r(k+a_m)/\alpha_m)} z^{-k/\alpha_m},
\ee
with the prime indicating the omission of the term corresponding to $r=m$ in the product. 
This expression in (\ref{e25}) consists of (at most) $p$ expansions each with the leading behaviour 
$z^{-a_m/\alpha_m}$ ($1\leq m\leq p$).
When the parameters $\alpha_r$ and $a_r$ are such that some of the poles 
are of higher order, the expansion (\ref{e25a}) is invalid and the residues must 
then be evaluated according to the multiplicity of the poles concerned; this will lead to terms involving $\log\,z$ in the algebraic expansion.

The three main expansion theorems are as follows. Throughout we let $\epsilon$ denote an arbitrarily 
small positive quantity.
\newtheorem{theorem}{Theorem}
\begin{theorem}$\!\!\!.$
If $0<\kappa<2$, then
\bee\label{e21}{}_p\Psi_q(z)\sim\left\{\begin{array}{lll}E_{p,q}(z)+H_{p,q}(ze^{\mp\pi i}) & 
\mbox{in} & |\arg\,z|\leq \f{1}{2}\pi\kappa \\
\\H_{p,q}(ze^{\mp\pi i}) & 
\mbox{in} & \fs\pi\kappa+\epsilon\leq|\arg\,z|\leq\pi\end{array} \right.
\ee
as $|z|\ra\infty$. The upper or lower sign 
in $H_{p,q}(ze^{\mp\pi i})$ is chosen according as $\arg\,z>0$ or $\arg\,z<0$, respectively.
\end{theorem}
It is seen that the $z$-plane is 
divided into two sectors, with a common vertex at $z=0$, by the rays  
$\arg\,z=\pm\f{1}{2}\pi\kappa$. In the sector $|\arg\,z|<\f{1}{2}\pi\kappa$, 
the asymptotic character of $_p\Psi_q(z)$ is exponentially large whereas in the complementary sector 
$\fs\pi\kappa<|\arg\,z|\leq\pi$, the dominant expansion of $_p\Psi_q(z)$ is algebraic in character. On the rays $\arg\,z=\pm\fs\pi\kappa$ the exponential expansion is oscillatory and is of a comparable magnitude to $H_{p,q}(ze^{\mp\pi i})$.
\begin{theorem}$\!\!\!.$
If $\kappa=2$ then
\bee\label{e22}
{}_p\Psi_q(z)\sim E_{p,q}(z)+E_{p,q}(ze^{\mp2\pi i})+H_{p,q}(ze^{\mp\pi i})
\ee
as $|z|\to\infty$ in the sector $|\arg\,z|\leq\pi$. The upper or lower signs are chosen according as $\arg\,z>0$ or $\arg\,z<0$, respectively.
\end{theorem}
The rays $\arg\,z=\pm\fs\pi\kappa$ now coincide with the negative real axis. It follows that ${}_p\Psi_q(z)$ is exponentially large in character as $|z|\to\infty$ except in the neighbourhood of the negative real axis, where the algebraic expansion becomes asymptotically significant.
\begin{theorem}$\!\!\!.$
When $\kappa>2$ we have\footnote{In \cite{W1}, the expansion was given in terms of the two dominant expansions only, viz. $E_{p,q}(z)$ and $E_{p,q}(ze^{\mp2\pi i})$, corresponding to $n=0$ and $n=\pm 1$ in (\ref{e23}).} 
\bee\label{e23}
{}_p\Psi_q(z)\sim \sum_{n=-N}^N E_{p,q}(ze^{2\pi in})+H_{p,q}(ze^{\mp\pi iz})
\ee
as $|z|\to\infty$ in the sector $|\arg\,z|\leq\pi$. The integer $N$ is chosen such that it is the smallest integer satisfying $2N+1>\fs\kappa$ and the upper or lower is chosen according as $\arg\,z>0$ or $\arg\,z<0$, respectively.

In this case the asymptotic behaviour of ${}_p\Psi_q(z)$ is exponentially large for all values of $\arg\,z$ and, consequently, the algebraic expansion may be neglected. The sums $E_{p,q}(ze^{2\pi in})$ are exponentially large (or oscillatory) as $|z|\to\infty$ for values of $\arg\,z$ satisfying $|\arg\,z+2\pi n|\leq\fs\pi\kappa$. 
\end{theorem}

The division of the $z$-plane into regions where ${}_p\Psi_q(z)$ possesses exponentially large or algebraic behaviour for large $|z|$ is illustrated in Fig.~1.
When $0<\kappa<2$, the exponential expansion $E_{p,q}(z)$ is still present in the sectors $\fs\pi\kappa<|\arg\,z|<\min\{\pi,\pi\kappa\}$, where it is subdominant. The rays $\arg\,z=\pm\pi\kappa$ ($0<\kappa<1$), where $E_{p,q}(z)$ is {\it maximally} subdominant with respect to $H_{p,q}(ze^{\mp\pi i})$, are called Stokes lines.\footnote{The positive real axis $\arg\,z=0$ is also a Stokes line where the algebraic expansion is maximally subdominant.} As these rays are crossed (in the sense of increasing $|\arg\,z|$) the exponential expansion switches off according to Berry's now familiar error-function smoothing law \cite{B}; see \cite{P10} for details. The rays $\arg\,z=\pm\fs\pi\kappa$, where $E_{p,q}(z)$ is oscillatory and comparable to $H_{p,q}(ze^{\mp\pi i})$, are called anti-Stokes lines.

\begin{figure}[th]
	\begin{center}
%{\small($a$)}
\includegraphics[width=0.35\textwidth]{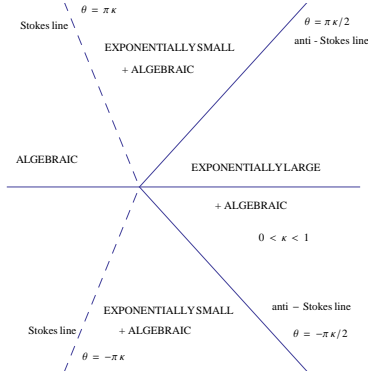}
%\qquad{\small($b$)}\includegraphics[width=0.35\textwidth]{expsWFnFig1A}
\caption{\small{The exponentially large and algebraic sectors associated with ${}_p\Psi_q(z)$ in the complex $z$-plane with $\theta=\arg\,z$ when $0<\kappa<1$. The Stokes and anti-Stokes lines are indicated.}}
	\end{center}
\end{figure}

In view of the above interpretation of the Stokes phenomenon a more precise version of Theorem 1 is as follows:
\begin{theorem}$\!\!\!.$
When $0<\kappa\leq 2$, then 
\bee\label{e24}
{}_p\Psi_q(z)\sim\left\{\begin{array}{lll}
E_{p,q}(z)+H_{p,q}(ze^{\mp\pi i}) & 
\mbox{in} & |\arg\,z|\leq \min\{\pi-\epsilon,\pi\kappa-\epsilon\} \\
%\\ \fs {\cal E}_{p,q}(z)+H_{p,q}^o(ze^{\mp\pi i}) & \mbox{on} & \arg\,z=\pm\pi\kappa\\ 
\\H_{p,q}(ze^{\mp\pi i}) & 
\mbox{in} & \pi\kappa+\epsilon\leq |\arg\,z|\leq\pi\ \ (0<\kappa<1)\\
\\E_{p,q}(z)+E_{p,q}(ze^{\mp2\pi i})+H_{p,q}(ze^{\mp\pi i}) & 
\mbox{in} & |\arg\,z|\leq\pi\ \ (1<\kappa\leq 2)\end{array} \right.
\ee
as $|z|\ra\infty$. The upper or lower signs 
are chosen according as $\arg\,z>0$ or $\arg\,z<0$, respectively. 
\end{theorem}
We omit the expansion {\it on} the Stokes lines $\arg\,z=\pm\pi\kappa$; the details in the case $p=1$, $q\geq0$ are discussed in \cite{P14}.
The expansions in (\ref{e24}a) and (\ref{e21}a) were given by Wright \cite{W1, W2} in the sector 
$|\arg\,z|\leq\min\{\pi, \f{3}{2}\pi\kappa-\epsilon\}$ as he did not take into account the Stokes phenomenon.
Since $E_{p,q}(z)$ is exponentially small in $\fs\pi\kappa<|\arg\,z|\leq\pi$, then in the sense of Poincar\'e, the expansion $E_{p,q}(z)$ can be neglected and there is no inconsistency between Theorems 1 and 4.
Similarly, $E_{p,q}(ze^{-2\pi i})$ is exponentially small compared to $E_{p,q}(z)$ in $0\leq\arg\,z<\pi$ and there is no inconsistency between the expansions in (\ref{e21}a) and (\ref{e24}c) when $1<\kappa<2$. However, in the vicinity of $\arg\,z=\pi$, these last two expansions are of comparable magnitude and, for real parameters, they combine to generate a real result on this ray. A similar remark applies to $E_{p,q}(ze^{2\pi i})$ in $-\pi<\arg\,z\leq 0$. 

The following theorem was given by Braaksma \cite[p.~331]{Br}.
\begin{theorem}$\!\!\!.$
If $p=0$, so that $g(s)$ has no poles and $\kappa>1$, then $H_{0,q}(z)\equiv 0$. When $1<\kappa<2$, we have the expansion 
\bee\label{e26}
{}_0\Psi_q(z)\sim E_{0,q}(z)+E_{0,q}(ze^{\mp2\pi i})
\ee
as $|z|\to\infty$ in the sector $|\arg\,z|\leq\pi$ The upper or lower sign is chosen according as $\arg\,z>0$ or $\arg\,z<0$, respectively. The dominant expansion ${}_0\Psi_q(z)\sim E_{p,q}(z)$ holds in the reduced sector $|\arg\,z|\leq\pi-\epsilon$.
\end{theorem}
It can be seen that (\ref{e26}) agrees with (\ref{e24}c) when $H_{p,q}(z)\equiv 0$. Braaksma gave the result (\ref{e26}) valid in a sector straddling the negative real axis given by $\pi-\delta\leq\arg\,z\leq\pi+\delta$, where $0<\delta<\fs\pi(1-\fs\kappa)$.

It is our purpose here to examine Theorems 4 and 5 in more detail by means of a series of examples. We carry out a numerical investigation to show that (\ref{e24}c) is valid when $1<\kappa<2$ and, when $0<\kappa<1$, that the exponential expansion $E_{p,q}(z)$ in Theorem 4 switches off (as $|\arg\,z|$ increases) across the Stokes lines $\arg\,z=\pm\pi\kappa$, where $E_{p,q}(z)$ is maximally subdoiminant with respect to $H_{p,q}(ze^{\mp\pi i})$. Similarly in Theorem 5, we show that when $1<\kappa<2$ the expansions $E_{p,q}(ze^{\mp2\pi i})$ switch off across the Stokes lines $\arg\,z=\pm\pi(1-\fs\kappa)$, where they are maximally subdominant with respect to $E_{p,q}(z)$. Thus, although the expansions in (\ref{e24}a) and (\ref{e26}) are valid asymptotic descriptions, more accurate evaluation will result from taking into account the Stokes phenomenon as the above-mentioned rays are crossed.

\vspace{0.6cm}

\begin{center}
{\bf 3. Numerical examples}
\end{center}
\setcounter{section}{3}
\setcounter{equation}{0}
\renewcommand{\theequation}{\arabic{section}.\arabic{equation}}

\noindent{\bf Example 3.1}\ \ \ Our first example is the Mittag-Leffler function ${\cal E}_{a,b}(z)$ defined by
\[{\cal E}_{a,b}(z)=\sum_{n=0}^\infty\frac{z^n}{\g(an+b)},\]
where we consider $a>0$. This corresponds to a case of ${}_1\Psi_1(z)$ with the parameters $\kappa=a$, $h=a^{-a}$, $\vartheta=1-b$ and $g(s)=\g(1+s)/\g(as+b)$. Then from (\ref{e22c})--(\ref{e22b}), we have $Z=z^{1/a}$, $A_0=1/a$ with
$A_j=0$ for $j\geq 1$. The exponential and algebraic expansions are from (\ref{e22c}), (\ref{e25}) and (\ref{e25a}) given by
\[E_{1,1}(z)=\frac{1}{a} z^{(1-b)/a} \exp [z^{1/a}],\qquad H_{1,1}(ze^{\mp\pi i})=-\sum_{k=1}^\infty \frac{z^{-k}}{\g(b-ak)}.\]
Then, from Theorems 2, 3 and 4 we obtain the following asymptotic expansions\footnote{When $a=1$ we have ${\cal E}_{1,b}(z)=z^{1-b}e^z P(b-1,z)$, where $P(\alpha,z)=\gamma(\alpha,z)/\g(\alpha)$ is the normalised incomplete gamma function. It then follows from \cite[(8.2.5), (8.11.2)]{DLMF} that the expansion of ${\cal E}_{1,b}(z)$ is given by (\ref{e31}a) as $|z|\to\infty$ in $|\arg\,z|\leq\pi$.} as $|z|\to\infty$.

(i)\ When $0<a<1$
\bee\label{e31}
{\cal E}_{a,b}(z)\sim \left\{\begin{array}{ll} \dfrac{1}{a}z^{(1-b)/a} \exp [z^{1/a}]-\sum_{k=1}^\infty \dfrac{z^{-k}}{\g(b-ak)}& (|\arg\,z|\leq\pi a-\epsilon)\\
-{\displaystyle\sum_{k=1}^\infty \frac{z^{-k}}{\g(b-ak)}} & (\pi a+\epsilon\leq\arg\,z\leq\pi);\end{array}\right.
\ee

(ii)\ when $1<a<2$
\bee\label{e32}
{\cal E}_{a,b}(z)\sim \left\{\begin{array}{ll} \dfrac{1}{a}z^{(1-b)/a} \exp [z^{1/a}]-\sum_{k=1}^\infty \dfrac{z^{-k}}{\g(b-ak)}& (|\arg\,z|\leq\pi-\epsilon)\\
\dfrac{1}{a}z^{(1-b)/a} \exp [z^{1/a}]+\dfrac{1}{a}(ze^{\mp2\pi i})^{(1-b)/a} \exp [(ze^{\mp2\pi i})^{1/a}] & \\
-{\displaystyle\sum_{k=1}^\infty \frac{z^{-k}}{\g(b-ak)}} & (|\arg\,z|\leq\pi);\end{array}\right.
\ee

(iii)\ when $a=2$
\[{\cal E}_{a,b}(z)\sim\dfrac{1}{a}z^{(1-b)/a} \exp [z^{1/a}]+\dfrac{1}{a}(ze^{\mp2\pi i})^{(1-b)/a} \exp [(ze^{\mp2\pi i})^{1/a}]\]
\bee\label{e33}
-{\displaystyle\sum_{k=1}^\infty \frac{z^{-k}}{\g(b-ak)}} \qquad\qquad (|\arg\,z|\leq\pi);\hspace{2cm}
\ee

(iv)\ when $a>2$
\bee\label{e34}
{\cal E}_{a,b}(z)\sim \frac{1}{a}\sum_{n=-N}^N (ze^{2\pi in})^{(1-b)/a} \exp [z^{1/a}e^{2\pi in/a}]-\sum_{k=1}^\infty \frac{z^{-k}}{\g(b-ak)} \qquad (|\arg\,z|\leq\pi),
\ee
where $N$ is the smallest integer\footnote{The more refined treatment of ${\cal E}_{a,1}(z)$ discussed in \cite[Section 5.1.4]{PK} has the integer $N$ satisfying $N<\fs a<N+1$. The additional exponential expansions present in (\ref{e34}) with this choice of $N$ are, however, exponentially small for $|\arg\,z|\leq\pi$.} satisfying $2N+1>\fs a$.
The upper or lower signs are taken according as $\arg\,z>0$ or $\arg\,z<0$, respectively.

When $0<a<1$, it is established in \cite{P02} (see also \cite{WZ02}) that the exponential term $a^{-1} \exp [z^{1/a}]$ in (\ref{e31}a)is multiplied by the approximate factor involving the error function
\[\frac{1}{2}+\frac{1}{2} \mbox{erf} \bl[\frac{(\pi a\mp\theta)}{a}\,\sqrt{\frac{|z|}{2}}\br]\]
as $|z|\to\infty$ in the neighbourhood of the Stokes lines $\theta=\arg\,z=\pm\pi a$, respectively, where it is maximally subdominant. This shows that
the above exponential term indeed switches off in the familiar manner \cite{B} as one crosses the Stokes lines in the sense of increasing $|\theta|$ and that consequently the expansion in (\ref{e31}a) is valid in $|\arg\,z|\leq \pi a-\epsilon$.

On the negative real axis we put $z=-x$, with $x>0$.
From (\ref{e32}), we have when $1<a<2$
\[{\cal E}_{a,b}(-x)\sim\frac{1}{a}(xe^{\pi i})^{(1-b)/a} \exp [x^{1/a}e^{\pi i/a}]+\frac{1}{a}(xe^{-\pi i})^{(1-b)/a} \exp [(x^{1/a}e^{-\pi i/a}]\]
\[\hspace{4cm}-\sum_{k=1}^\infty \frac{(-x)^{-k}}{\g(b-ak)}\]
\bee\label{e35}
=F_{a,b}(x)-\sum_{k=1}^\infty \frac{(-x)^{-k}}{\g(b-ak)},
\ee
as $x\to+\infty$, where
\bee\label{e35a}
F_{a,b}(x)=\frac{2}{a} x^{(1-b)/a} \exp \bl[x^{1/a} \cos \frac{\pi}{a}\br]\,\cos \bl[x^{1/a} \sin \frac{\pi}{a}+\pi \frac{(1-b)}{a}\br].
\ee
 The presence of the additional exponential expansion $E_{1,1}(ze^{\mp2\pi i})$ in (\ref{e32}) is seen to be essential in order to obtain a real result\footnote{We remark that the result (\ref{e35}) can also be deduced by use of the identity
${\cal E}_{a,b}(-x)=2{\cal E}_{2a,b}(x^2)-{\cal E}_{a,b}(x)$ combined with the expansions of ${\cal E}_{a,b}(z)$ for $z\to+\infty$.}
(when $b$ is real) on the negative $z$-axis.

\vspace{0.4cm}

\noindent{\bf Example 3.2}\ \ \ Our second example is the function
\bee\label{e37}
F_1(z)=\sum_{n=0}^\infty \frac{\g(\fs n+a)}{\g(n+b)}\,\frac{z^n}{n!}\qquad (\kappa=\f{3}{2}),
\ee
where $a$ and $b$ are finite parameters, which corresponds to a case of ${}_1\Psi_1(z)$. The exponential expansion is
\[E_{1,1}(z)=Z^\vartheta e^Z \sum_{j=0}^\infty A_j Z^{-j}, \qquad Z=\f{3}{2}(hz)^{2/3},\]
where, from (\ref{e22b}),
\[A_0=(\f{2}{3})^{\vartheta+1/2} (\fs)^{a-1/2}\]
and $\vartheta=a-b$, $h=2^{-1/2}$. An algorithm for the computation of the normalised coefficients $c_j=A_j/A_0$ is described in the appendix. In our computations we have employed $0\leq j\leq 40$; the first ten coefficients $c_j$ for $F_1(z)$ are listed in Table 1 for the particular case $a=\f{1}{4}$ and $b=\f{3}{4}$. From (\ref{e25}), the algebraic expansion is
\[H_{1,1}(ze^{\mp\pi i})=2\sum_{k=0}^\infty\frac{(-)^k \g(2k+2a)}{k! \g(b-2a-2k)}\,(ze^{\mp\pi i})^{-2k-2a}.\]

\begin{table}[h]
\begin{center}
\begin{tabular}{l|l||l|l}
\mcol{1}{c|}{$j$} & \mcol{1}{c||}{$c_j$} & \mcol{1}{c|}{$j$} & \mcol{1}{c}{$c_j$}\\
[0.1cm]\hline
&&& \\[-0.2cm]
1 & $\f{61}{192}$ & 2 & $\f{23161}{73728}$\\
&&& \\[-0.2cm]
3 & $\f{22783285}{42467328}$ & 4 & $\f{44604509425}{32614907904}$\\
&&& \\[-0.2cm]
5 & $\f{30375638199305}{6262062317568}$& 6 & $\f{162721816250787605}{7213895789838336}$ \\
&&& \\[-0.2cm]
7 & $\f{180090830597703240215}{1385067991648960512}$ & 8 & $\f{1889199431389108590226475}{2127464435172803346432}$\\
&&& \\[-0.2cm]
9 & $\f{25599447910539396612172829375}{3676258543978604182634496}$ & 10 & $\f{86726322340809175676137010099575}{1411683280887784006131646464}$\\
[0.2cm]\hline
\end{tabular}
\end{center}
\caption{\footnotesize{The normalised coefficients $c_j$ for $1\leq j\leq10$ (with $c_0=1$) for the sum (\ref{e37}) when $a=\f{1}{4}$ and $b=\f{3}{4}$.}}
\end{table}

It is clearly sufficient for real parameters to consider values of $z$ satisfying $0\leq\arg\,z\leq\pi$ and this we do throughout this section. From Theorem 4, we obtain
\[F_1(z)=E_{1,1}(z)+E_{1,1}(ze^{-2\pi i})+H_{1,1}(ze^{-\pi i})\]
as $|z|\to\infty$ in $0\leq\arg\,z\leq\pi$, from which we see that $F_1(z)$ is exponentially large in the sector $|\arg\,z|<3\pi/4$. We have computed $F_1(z)$ for a value of $|z|$ and varying $\theta=\arg\,z$ in the range $0.7\pi\leq\theta\leq\pi$. In Table 2 we show the absolute values of
\begin{table}[h]
\begin{center}
\begin{tabular}{l|ll}
\mcol{1}{c|}{$\theta/\pi$} & \mcol{1}{c}{$|{\cal R}_1(z)|$} 
& \mcol{1}{c}{$|E_{1,1}(ze^{-2\pi i})|$}\\
[.1cm]\hline
&&\\[-0.25cm]
1.00  & $6.283513\times 10^{-7}$  & $6.283515\times 10^{-7}$\\
0.95  & $6.605074\times 10^{-8}$  & $6.605098\times 10^{-8}$ \\
0.90  & $8.190985\times 10^{-9}$  & $8.190854\times 10^{-9}$\\
0.85  & $1.226317\times 10^{-9}$  & $1.225981\times 10^{-9}$\\
0.80  & $2.263874\times 10^{-10}$ & $2.261409\times 10^{-10}$\\
0.75  & $5.240704\times 10^{-11}$ & $5.236698\times 10^{-11}$\\
0.70  & $1.573812\times 10^{-11}$ & $1.546959\times 10^{-11}$\\
[.2cm]\hline
\end{tabular}
\end{center}
\caption{\footnotesize{Values of the absolute error ${\cal R}_1(z)$ in the computation of $F_1(z)$
using an optimal truncation of both $E_{1,1}(z)$ and $H_{1,1}(ze^{-\pi i})$ compared with $|E_{1,1}(ze^{-2\pi i})|$ as a function of $\theta$ for $z=100e^{i\theta}$, 
$a=\f{1}{4}$ and $b=\f{3}{4}$.}}
\end{table}
\[{\cal R}_1(z)\equiv F_1(z)-E_{1,1}^{opt}(z)-H_{1,1}^{opt}(ze^{-\pi i})\]
compared with $|E_{1,1}(ze^{-2\pi i})|$ (which was computed for $0\leq j\leq 5$), where the superscript `opt' denotes that both the asymptotic sums $E_{1,1}(z)$ and $H_{1,1}(ze^{-\pi i})$ are truncated at their respective optimal truncation points. The results clearly confirm that (i) the exponential expansion $E_{1,1}(z)$ is present in the algebraic sector $\f{3}{4}\pi<\arg\,z\leq\pi$ and (ii) 
the subdominant expansion $E_{1,1}(ze^{-2\pi i})$ is present in (at least) the sector $0.7\pi\leq\theta\leq\pi$.
It was not possible to penetrate very far into the exponentially large sector $|\arg\,z|<\f{3}{4}\pi$, since the error in the computation of $E_{1,1}(z)$ --- even at optimal truncation --- swamps the algebraic and subdominant exponential expansions. Such a computation would require a hyperasymptotic evaluation of the dominant expansion on the lines of that described for the generalised Bessel function ${}_0\Psi_1(z)$ in Wong and Zhao \cite{WZ}.
\vspace{0.4cm}

\noindent{\bf Example 3.3}\ \ \ Consider the function
\[F_2(z)=\sum_{n=0}^\infty \frac{\g(\f{2}{3}n+a\,)z^n}{\g(\f{1}{3}n+b)\,n!}\qquad (\kappa=\f{2}{3}).\]
According to Theorem 4, the expansion of $F_2(z)$ for large $|z|$ is
\[F_2(z)\sim E_{1,1}(z)+H_{1,1}(ze^{-\pi i})\qquad (0\leq\arg\,z\leq \f{2}{3}\pi-\epsilon).\] 
The algebraic expansion is, from (\ref{e25}), given by
\[H_{1,1}(ze^{-\pi i})=\frac{3}{2}\sum_{k=0}^\infty\frac{(-)^k\g(\f{3}{2}k+\f{3}{2}a)}
{k!\,\g(b-\fs a-\fs k)}(ze^{-\pi i})^{-3(k+a)/2}\]
and the exponential expansion $E_{1,1}(z)$ is obtained from (\ref{e22c}) with the parameters $\vartheta=a-b$, $h=(\f{2}{3})^{\frac{2}{3}}(\f{1}{3})^{-\frac{1}{3}}$
and $A_0=\kappa^{-\fr-\vartheta}(\f{2}{3})^{a-\fr}(\f{1}{3})^{\fr-b}$. The coefficients $A_j$ are obtained as indicated in Example 3.2.

The function $F_2(z)$ is exponentially large in the sector $|\arg\,z|<\f{1}{3}\pi$, whereas in the
sector $\f{1}{3}\pi<\arg\,z\leq\pi$ the algebraic expansion $H_{1,1}(ze^{-\pi i})$ is dominant.
The expansion $E_{1,1}(z)$ is maximally subdominant with respect to $H_{1,1}(ze^{-\pi i})$ on the ray $\arg\,z=\pi\kappa=\f{2}{3}\pi$.
Consequently, as $\arg\,z$ increases, the exponential expansion $E_{1,1}(z)$ should switch off across the
Stokes line $\arg\,z=\f{2}{3}\pi$, to leave the algebraic expansion $H_{1,1}(ze^{-\pi i})$ in the sector
$\f{2}{3}\pi<\arg\,z\leq\pi$. To demonstrate this, we define the Stokes multiplier $S(\theta)$ by
\[F_2(z)=H_{11}^{opt}(ze^{-\pi i})+A_0Z^\vartheta e^Z\,S(\theta).\]
In Table 3 we show the absolute values of ${\cal R}_2(z):=F_2(z)-H_{1,1}^{opt}(ze^{-\pi i})$ and of the leading term of $E_{1,1}(z)$ as a function of $\theta=\arg\,z$. We also show the values\footnote{The Stokes multiplier $S(\theta)$ has a small imaginary part that we do not show.} 
of Re($S(\theta)$) in the neighbourhood of the Stokes line $\arg\,z=\f{2}{3}\pi$ for the case
$z=10e^{i\theta}$ and $a=\f{1}{3}$, $b=\f{1}{4}$. It is seen that the Stokes multiplier has the value $\simeq 1$ when $\theta=\fs\pi$ (before the transition commences) and $\simeq 0$ when $\theta=\f{3}{4}\pi$ (after the transition is almost completed).
\begin{table}[h]
\begin{center}
\begin{tabular}{l|llc}
\multicolumn{1}{c|}{$\theta/\pi$} & \multicolumn{1}{c}{$|{\cal R}_2(z)|$} 
& \mcol{1}{c}{$|A_0Z^\vartheta e^Z|$} & \mcol{1}{c}{Re$(S(\theta))$}\\
[.1cm]\hline
&&\\[-0.25cm]
0.50  & $4.4964\times 10^{-8}$  & $4.4947\times 10^{-8}$  & 1.0000 \\
0.55  & $1.2980\times 10^{-9}$  & $1.3005\times 10^{-9}$  & 0.9981 \\
0.60  & $1.1196\times 10^{-10}$ & $1.1848\times 10^{-10}$ & 0.9450 \\
0.62  & $5.6361\times 10^{-11}$ & $6.4685\times 10^{-11}$ & 0.8713 \\
0.64  & $3.2641\times 10^{-11}$ & $4.3607\times 10^{-11}$ & 0.7485 \\
0.66  & $1.9737\times 10^{-11}$ & $3.6426\times 10^{-11}$ & 0.5418 \\
0.68  & $1.3545\times 10^{-11}$ & $3.7762\times 10^{-11}$ & 0.3600 \\
0.70  & $9.9952\times 10^{-12}$ & $4.8568\times 10^{-11}$ & 0.2058 \\
0.72  & $9.1973\times 10^{-12}$ & $7.7328\times 10^{-11}$ & 0.1189 \\
0.75  & $5.6314\times 10^{-12}$ & $2.2959\times 10^{-10}$ & 0.0237 \\
[.2cm]\hline
\end{tabular}
\end{center}
\caption{\footnotesize{Values of the absolute error in ${\cal R}_2(z)\equiv F_2(z)-H_{1,1}^{opt}(ze^{-\pi i})$ in the computation of $F_2(z)$
using an optimal truncation of the algebraic expansion compared with the leading term of $|E_{1,1}(z)|$ as a function of $\theta$ for $z=10e^{i\theta}$, 
$a=\f{1}{3}$ and $b=\f{1}{4}$. The final column shows the real part of the computed Stokes multiplier $S(\theta)$ for transition across the ray $\arg\,z=\f{2}{3}\pi$.}}
\end{table}

\vspace{0.4cm}

\noindent{\bf Example 3.4}\ \ \ Our final example is the function of the type ${}_0\Psi_2(z)$ given by
\bee
F_3(z)=\sum_{n=0}^\infty\frac{z^n}{n! \g(cn+a) \g(cn+b)}\qquad (\kappa=1+2c),
\ee
where $0<c\leq\fs$. Since $p=0$, the algebraic expansion $H_{0,2}(z)\equiv 0$. From Theorem 5 we obtain the asymptotic expansion
\[F_3(z)\sim E_{0,2}(z)+E_{0,2}(ze^{\mp2\pi i})\qquad (|\arg\,z|\leq\pi),\]
where the associated parameters are $\vartheta=1-a-b$, $h=c^{-2c}$ and 
\[A_0=\frac{c^\vartheta \kappa^{-\vartheta-1/2}}{2\pi}~.\]
The function $F_3(z)$ is exponentially large in the sector $|\arg\,z|<\fs\pi(1+2c)$. The other expansion $E_{0,2}(ze^{-2\pi i})$ is subdominant in the upper half-plane but combines with $E_{0,2}(z)$ on the negative real axis to produce (for real $a$ and $b$) a real expansion.

Since the exponential factors associated with $E_{0,2}(z)$ and $E_{0,2}(ze^{-2\pi i})$ are $\exp [|Z|e^{i\theta/\kappa}]$ and $\exp [|Z|e^{i(\theta-2\pi)/\kappa}]$, where $\theta=\arg\,z$ and we recall that $Z$ is defined in (\ref{e22c}), the greatest difference between these factors occurs when
\[\sin \bl(\frac{\theta}{\kappa}\br)=\sin \bl(\frac{\theta-2\pi}{\kappa}\br);\]
that is, when $\theta=\fs\pi(2-\kappa)$. Consequently, as  $\arg\,z$ increases in the upper half-plane, we expect that the expansion $E_{0,2}(ze^{-2\pi i})$ should switch on across the Stokes line $\arg\,z=\fs\pi(2-\kappa)$; similar considerations apply to $E_{0,2}(ze^{2\pi i})$ and the Stokes line $\arg\,z=-\fs\pi(2-\kappa)$ in the lower half-plane.

To demonstrate the correctness of this claim, we choose $c=\f{1}{10}$ (so that $\kappa=\f{6}{5}$) and $a=\f{1}{4}$, $b=\f{3}{4}$. The function $F_3(z)$ is therefore exponentially large in the sector $|\arg\,z|<\f{3}{5}\pi$ and the Stokes line in the upper half-plane is $\arg\,z=\f{2}{5}\pi$. We have chosen $a-b$ to have a half-integer value for a very specific reason. The more detailed treatment in \cite{P10} shows that there is a {\it third} ({\it subdominant\/}) {\it exponential series} present in the expansion of $F_3(z)$ given by
\[2\cos \pi(a-b)\,X^\vartheta e^{-X} \sum_{j=0}^\infty A_j (-X)^{-j},\qquad X=\kappa(hze^{-\pi i})^{1/\kappa}.\]
Our present choice of $a$ and $b$ therefore eliminates this third expansion and enables us to deal with a case comprising only two exponential expansions.

In Table 4, we show for $|z|=20$ and varying $\theta=\arg\,z$ the values of $|F_3(z)-E_{0,2}^{opt}(z)|$ and $|E_{0,2}(ze^{-2\pi i})|$ together with the real part of the Stokes multiplier $S(\theta)$ defined by
\[F_3(z)=E_{0,2}^{opt}(z)+A_0 (Ze^{-2\pi i/\kappa})^\vartheta\,\exp [Ze^{-2\pi i/\kappa}]\,S(\theta).\] 
The results clearly demonstrate the switching-on of the subdominant expansion $E_{0,2}(ze^{-2\pi i})$ across the Stokes line $\arg\,z=\f{2}{5}\pi$ as $\arg\,z$ increases in the upper half-plane. 
\begin{table}[h]
\begin{center}
\begin{tabular}{l|ccc}
\multicolumn{1}{c|}{$\theta/\pi$} & \multicolumn{1}{c}{$|{\cal R}_3(z)|$} 
& \multicolumn{1}{c}{$|E_{0,2}(ze^{-2\pi i})|$} & \multicolumn{1}{c}{Re$(S(\theta))$}\\
[.1cm]\hline
&&\\[-0.25cm]
0.20  & $7.231938\times 10^{-4}$ & $1.452127\times 10^{-1}$ & 0.0020 \\
0.25  & $2.204854\times 10^{-4}$ & $8.898995\times 10^{-3}$ & 0.0184 \\
0.30  & $5.082653\times 10^{-5}$ & $5.720603\times 10^{-4}$ & 0.0797 \\
0.35  & $9.416276\times 10^{-6}$ & $4.042959\times 10^{-5}$ & 0.2230 \\
0.40  & $1.502207\times 10^{-6}$ & $3.287009\times 10^{-6}$ & 0.4477 \\
0.45  & $2.239289\times 10^{-7}$ & $3.209167\times 10^{-7}$ & 0.6893 \\
0.50  & $3.430029\times 10^{-8}$ & $3.915246\times 10^{-8}$ & 0.8679 \\
0.55  & $5.977355\times 10^{-9}$ & $6.187722\times 10^{-9}$ & 0.9575 \\
0.60  & $1.301304\times 10^{-9}$ & $1.307416\times 10^{-9}$ & 0.9862 \\
1.00  & $1.307416\times 10^{-9}$ & $1.307416\times 10^{-9}$ & 0.9908 \\
[.1cm]\hline
\end{tabular}
\end{center}
\caption{\footnotesize{Values of the absolute error in ${\cal R}_3(z)\equiv F_3(z)-E_{0,2}^{opt}(z)$ in the computation of $F_3(z)$
using an optimal truncation of $E_{0,2}(z)$ compared with $|E_{0,2}(ze^{-2\pi i})|$ as a function of $\theta$ for $z=20e^{i\theta}$, 
$a=\f{1}{4}$ and $b=\f{3}{4}$. The final column shows the real part of the computed Stokes multiplier $S(\theta)$ for transition across the ray $\arg\,z=\f{2}{5}\pi$.}}
\end{table}

\vspace{0.6cm}

\begin{center}
{\bf Appendix: \ An algorithm for the computation of the coefficients $c_j=A_j/A_0$}
\end{center}
\setcounter{section}{1}
\setcounter{equation}{0}
\renewcommand{\theequation}{\Alph{section}.\arabic{equation}}
We describe an algorithm for the computation of the normalised coefficients $c_j=A_j/A_0$ appearing in the exponential expansion $E_{p,q}(z)$ in (\ref{e22c}). Methods of computing these coefficients by recursion in the case $\alpha_r=\beta_r=1$ have been given by Riney \cite{R} and Wright \cite{W3}; see \cite[Section 2.2.2]{PK} for details. Here we describe an algebraic method for arbitrary $\alpha_r>0$ and $\beta_r>0$.

The inverse factorial expansion (\ref{e22a}) can be re-written as
\begin{equation}\label{a1}
\frac{g(s)\Gamma(\kappa s+\vartheta')}{\Gamma(1+s)}
=\kappa A_0(h\kappa^\kappa)^{s}\bl\{\sum_{j=0}^{M-1}\frac{c_j}{(\kappa s+\vartheta')_j}+\frac{O(1)}{(\kappa s+\vartheta')_M}\br\}
\end{equation}
for $|s|\to\infty$ uniformly in $|\arg\,s|\leq\pi-\epsilon$, where $g(s)$ is defined in (\ref{e11a}) with $n$ replaced by $s$. Introduction of the scaled gamma function $\g^*(z)=\g(z) (2\pi)^{-\fr}e^z z^{\fr-z}$ leads to the representation
\[\g(\alpha s+a)= (2\pi)^\fr e^{-\alpha s} (\alpha s)^{\alpha s+a-\fr} \,{\bf e}(\alpha s; a)\g^*(\alpha s+a),\]
where
\[{\bf e}(\alpha s; a):= e^{-a}\bl(1+\frac{a}{\alpha s}\br)^{\alpha s+a-\fr}=\exp\,\left[(\alpha s+a-\fs) \log\,\left(1+\frac{a}{\alpha s}\right)-a\right].\]

Then, after some routine algebra we find that the left-hand side of (\ref{a1}) can be written as
\bee\label{a2}
\frac{g(s) \g(\kappa s+\vartheta')}{\g(1+s)}=\kappa A_0(h\kappa^\kappa)^{s}\,R_{p,q}(s)\,\Upsilon_{p,q}(s),
\ee
where
\[\Upsilon_{p,q}(s):=\frac{\prod_{r=1}^p\g^*(\alpha_rs+a_r)}{\prod_{r=1}^q\g^*(\beta_rs+b_r)}\,\frac{\g^*(\kappa s+\vartheta')}{\g^*(1+s)},\qquad R_{p,q}(s):=\frac{\prod_{r=1}^p e(\alpha_rs;a_r)}{\prod_{r=1}^q e(\beta_rs;b_r)}\,\frac{e(\kappa s;\vartheta')}{e(s;1)}.\]
Substitution of (\ref{a2}) in (\ref{a1}) then yields the inverse factorial expansion in the form 
\bee\label{a3}
R_{p,q}(s)\,\Upsilon_{p,q}(s)=\sum_{j=0}^{M-1}\frac{c_j}{(\kappa s+\vartheta')_j}+\frac{O(1)}{(\kappa s+\vartheta')_M}
\ee
as $|s|\to\infty$ in $|\arg\,s|\leq\pi-\epsilon$.

We now expand $R_{p,q}(s)$ and $\Upsilon_{p,q}(s)$ for $s\to+\infty$ making use of the well-known expansion (see, for example, \cite[p.~71]{PK})
\[\g^*(z)\sim\sum_{k=0}^\infty(-)^k\gamma_kz^{-k}\qquad(|z|\ra\infty;\ |\arg\,z|\leq\pi-\epsilon),\]
where $\gamma_k$ are the Stirling coefficients, with 
\[\gamma_0=1,\quad \gamma_1=-\f{1}{12},\quad \gamma_2=\f{1}{288},\quad  \gamma_3=\f{139}{51840},
\quad \gamma_4=-\f{571}{2488320}, \ldots\ .\]
Then we find
\[\g^*(\alpha s+a)=1-\frac{\gamma_1}{\alpha s}+O(s^{-2}),\qquad e(\alpha s;a)=1+\frac{a(a-1)}{2\alpha s}+O(s^{-2}),\]
whence
\[R_{p,q}(s)=1+\frac{{\cal A}}{2s}+O(s^{-2}),\qquad
\Upsilon_{p,q}(s)=1+\frac{{\cal B}}{12s}+O(s^{-2}),\]
where we have defined the quantities ${\cal A}$ and ${\cal B}$ by
\[{\cal A}=\sum_{r=1}^p \frac{a_r(a_r-1)}{\alpha_r}-\sum_{r=1}^q\frac{b_r(b_r-1)}{\beta_r}-\frac{\vartheta}{\kappa} (1-\vartheta),\qquad
{\cal B}=\sum_{r=1}^p\frac{1}{\alpha_r}-\sum_{r=1}^q\frac{1}{\beta_r}+\frac{1}{\kappa}-1.\]
Upon equating coefficients of $s^{-1}$ in (\ref{a3}) we then obtain
\bee\label{a4}
c_1=\fs\kappa({\cal A}+\f{1}{6} {\cal B}).
\ee

The higher coefficients are obtained by continuation of this expansion process in inverse powers of $s$. We write the product on the left-hand side of (\ref{a3}) as an expansion in inverse powers of $\kappa s$ in the form
\bee\label{a5}
R_{p,q}(s) \Upsilon_{p,q}(s)=1+\sum_{j=1}^{M-1} \frac{C_j}{(\kappa s)^j}+O(s^{-M})
\ee
as $s\to+\infty$, where the coefficients $C_j$ are determined with the aid of {\it Mathematica}.
From the expansion of the ratio of two gamma functions in \cite[(5.11.13)]{DLMF} we obtain
\[\frac{1}{(\kappa s+\vartheta')_j}=\frac{1}{(\kappa s)^{j}} \bl\{\sum_{j=0}^{M-1}\frac{(-)^k(j)_k}{(\kappa s)^k k!}\,B_k^{(1-j)}(\vartheta')+O(s^{-M})\br\},\]
where $B_k^{(\sigma)}(x)$ are the generalised Bernoulli polynomials defined by
\[\bl(\frac{t}{e^t-1}\br)^{\sigma} e^{xt}=\sum_{k=0}^\infty \frac{B_k^{(\sigma)}(x)}{k!}\,t^k \qquad (|t|<2\pi).\]
Here we have $\sigma=1-j\leq 0$ and $B_0^{(\sigma)}(x)=1$. 

Then the right-hand side of (\ref{a3}) as $s\to+\infty$ becomes
\[1+\sum_{j=1}^{M-1}\frac{c_j}{(\kappa s+\vartheta')_j}+O(s^{-M})=1+\sum_{j=1}^{M-1}\frac{c_j}{(\kappa s)^j} \sum_{k=0}^{M-1} \frac{(-)^k (j)_k}{(\kappa s)^k k!}\,B_k^{(1-j)}(\vartheta')+O(s^{-M})\]
\bee\label{a6}
=1+\sum_{j=1}^{M-1} \frac{D_j}{(\kappa s)^j}+O(s^{-M})
\ee
with
\[D_j=\sum_{k=0}^{j-1} (-)^k\bl(\!\!\!\begin{array}{c}j-1\\k\end{array}\!\!\!\br) c_{j-k}\,B_k^{(k-j+1)}(\vartheta'),\]
where we have made the change in index $j+k\to j$ and used `triangular' summation (see \cite[p.~58]{S}).
Substituting (\ref{a5}) and (\ref{a6}) into (\ref{a3}) and equating the coefficients of like powers of $\kappa s$, we then find $C_j=D_j$ for $1\leq j\leq M-1$,
whence
\[c_j=C_j-\sum_{k=1}^{j-1}(-)^k\bl(\!\!\!\begin{array}{c}j-1\\k\end{array}\!\!\!\br) c_{j-k}\,B_k^{(k-j+1)}(\vartheta').\]

Thus we find 
\begin{eqnarray*}
c_1&=&C_1,\\
c_2&=&C_2-c_1 B_1^{(0)}(\vartheta'),\\
c_3&=&C_3-2c_2 B_1^{(-1)}(\vartheta')+c_1 B_2^{(0)}(\vartheta'), \ldots
\end{eqnarray*}
and so on, from which the coefficients $c_j$ can be obtained recursively.
With the aid of {\it Mathematica} this procedure is found to work well in specific cases when the various parameters have numerical values, where up to a maximum of 100 coefficients have been so calculated.


\vspace{0.6cm}

\begin{thebibliography}{99}
\footnotesize{
%\bibitem{} M. Abramowitz and I. A. Stegun, Handbook of Mathematical Functions, Dover, New York, 1965.

\bibitem{B} M.V. Berry, Uniform asymptotic smoothing of Stokes's discontinuities, Proc. Roy. Soc. London {\bf A422} (1989) 7--21.

\bibitem{Br} B.L.J. Braaksma, Asymptotic expansions and analytic continuations for a class of Barnes-integrals,
Compos. Math. {\bf 15} (1963) 239--341.

\bibitem{F}
W.B. Ford, {\it The Asymptotic Developments of Functions Defined by Maclaurin Series}, University of Michigan, Science Series, Vol. II, 1936.

\bibitem{N}
C.V. Newsom, On the character of certain entire functions in distant portions of the plane, Amer. J. Math. {\bf 60} (1938) 561--572. 

\bibitem{DLMF}
F.W.J. Olver, D.W. Lozier, R.F. Boisvert and C.W. Clark (eds.),    
{\it NIST Handbook of Mathematical Functions}, Cambridge University Press, Cambridge, 2010.

%\bibitem{P91} R.B. Paris, Smoothing of the Stokes phenomenon using Mellin-Barnes integrals,
%J. Comput. Appl. Math. {\bf 41} (1992) 117--133.

%\bibitem{P92} R.B. Paris, Smoothing of the Stokes phenomenon for high-order differential equations,
%Proc. Roy. Soc. London {\bf A436} (1992) 165--186.
\bibitem{P02} R.B. Paris, Exponential asymptotics of the Mittag-Leffler function, Proc. Roy. Soc. London {\bf 458A} (2002) 3041--3052.

\bibitem{P09} R.B. Paris, Some remarks on the theorems of Wright and Braaksma concerning the asymptotics of the generalised hypergeometric functions, Technical Report MS 09:01, Abertay University, 2009.

\bibitem{P10} R.B. Paris, Exponentially small expansions in the asymptotics of the Wright function, J. Comput. Appl. Math. {\bf 234} (2010) 488--504.

\bibitem{P14} R.B. Paris, Exponentially small expansions of the Wright function on the Stokes lines, Lithuanian Math. J. {\bf 54} (2014) 82--105..

\bibitem{PW} R.B. Paris and A.D. Wood, {\it Asymptotics of High Order Differential Equations}, Pitman Research 
Notes in Mathematics, {\bf 129}, Longman Scientific and Technical, Harlow, 1986.

\bibitem{PK} R.B. Paris and D. Kaminski,  {\it Asymptotics and Mellin-Barnes Integrals}, 
Cambridge University Press, Cambridge, 2001.

\bibitem{R}
T.D. Riney, On the coefficients in asymptotic factorial expansions, Proc. Amer. Math. Soc. {\bf 7} (1956) 245--249.

\bibitem{S} L.J. Slater, {\it Generalized Hypergeometric Functions}, Cambridge University Press, Cambridge, 1966.

\bibitem{WZ} R. Wong and Y.-Q. Zhao, Smoothing of Stokes's discontinuity for the generalized Bessel %function,
Proc. Roy. Soc. London {\bf A455} (1999) 1381--1400.

\bibitem{WZ02} R. Wong and Y.-Q. Zhao, Exponential asymptotics of the Mittag-Leffler function, Constr. Approx. {\bf 18} (2002) 355--385.

\bibitem{W1} E.M. Wright, The asymptotic expansion of the generalized hypergeometric function, J. Lond. Math. Soc. (Ser. 2) {\bf 10} (1935) 286--293.

\bibitem{W2} E.M. Wright, The asymptotic expansion of the generalized hypergeometric function,
Proc. Lond. Math. Soc. (Ser. 2) {\bf 46} (1940) 389--408.

\bibitem{W3} E.M. Wright, A recursion formula for the coefficients in an asymptotic expansion, Proc. Glasgow Math. Assoc. {\bf 4} (1958) 38--41. 
}
\end{thebibliography}
\end{document}